\documentclass[11pt]{article}
\usepackage{latexsym,amssymb,amsthm}
\usepackage{amsmath}
\def\s{\sigma}

\def\tbase{{\rm Base}}

\def\cT{{\cal T}}

\def\CC{\mathbb C}

\def\ZZ{\mathbb Z}

\def\PP{\mathbb P}

\def\QQ{\mathbb Q}

\def\fh{{\mathfrak h}}
\def\f4{{\mathfrak f}_4}
\def\fii{{\mathfrak i}}
\def\fsl{{\mathfrak {sl}}}

\def\fsp{{\mathfrak {sp}}}

\def\fso{{\mathfrak {so}}}
\def\fe{{\mathfrak e}}
\def\fa{{\mathfrak a}}
\def\fd{{\mathfrak d}}
\def\fg{{\mathfrak g}}

\def\ft{{\mathfrak t}}

\def\l{\lambda}
\def\a{\alpha}
\def\o{\omega}

\def\b{\beta}
\def\g{\gamma}

\def\m{\mu}
\def\r{\rho}

\def\ot{{\mathord{\otimes }\,}}
\def\op{{\mathord{\oplus }\,}}
\def\otc{{\mathord{\otimes\cdots\otimes}\;}} 

\def\ra{{\mathord{\;\rightarrow\;}}}
\def\we{\wedge}

 \def\BC{\mathbb C }

\def\BO{\mathbb O}\def\BS{\mathbb S }

\def\11{{\bf 1}}

\def\bcc#1{\mathbb C^{#1}}

\def\ff#1{\mathbb F\mathbb F^{#1}}
\def\fff#1#2{\mathbb F\mathbb F^{#1}_{#2}} \def\fg{\mathfrak g}
\def\fh{\mathfrak h}  \def\ft{\mathfrak t}

\def\hd{, ... ,}

\def\La#1{\Lambda^{#1}}

\def\ot{\!\otimes\!}

\def\pp#1{\mathbb P^{#1}}
\def\ppp{\mathbb P}
\def\qq#1#2#3{q^{#1}_{{#2} {#3}}}

\def\ra{\rightarrow}
\def\rr#1#2#3#4{r^{#1}_{{#2} {#3}{#4}}}
\def\sx{\sigma (X)}

\def\tbase{{\rm Base}\,}

\def\tdim{{\rm dim}\,}

\def\tcodim{{\rm codim}\,}

\def\tmod{{\rm mod }}

\def\th{\theta}
\def\TH{\Theta}

\def\up#1{{}^{({#1})}}

\def\ww{\wedge}
\def\ux#1{x^{#1}}

\def\xx#1{x^{#1}}

\begin{document}

\title{
Construction and classification of complex \linebreak 
simple Lie algebras \linebreak via projective geometry}
\author{ J.M. Landsberg and Laurent Manivel}
\date{  }
  \maketitle

\abstract{We construct the complex simple Lie algebras
using elementary algebraic geometry. We use our construction
to obtain a new proof of the classification of complex simple Lie
algebras that does not appeal to the classification of root systems.}
\endabstract

\newcommand{\norm}[1]{\lVert#1\rVert}

\theoremstyle{plain}
\newtheorem{theo}{Theorem}[section]
\newtheorem{coro}[theo]{Corollary}
\newtheorem{lemm}[theo]{Lemma}
\newtheorem{prop}[theo]{Proposition}

\theoremstyle{definition}
\newtheorem{defi}[theo]{Definition}

\theoremstyle{remark}
\newtheorem{rema}[theo]{Remark}
\newtheorem{exam}[theo]{Example}

\section{Overview}

We first present an algorithm that constructs the {\it minuscule
varieties} using elementary algebraic geometry. The minuscule
varieties are a preferred class of homogeneous varieties.
  They are essentially the homogeneous projective varieties that admit
an irreducible Hermitian symmetric metric; see below for the
precise definition. The algorithm proceeds iteratively by
  building a larger space $X\subset\pp
N$ from a smaller space $Y\subset\pp n$ via a rational map $\pp
n\dashrightarrow \pp N$, defined using the ideals of the secant  
  varieties of $Y$, beginning with $Y=\BC\pp 1$.
As a byproduct,  we obtain elementary constructions of all complex
simple Lie algebras  (except for $\frak e_8$  which has no minuscule
homogenous space) and their minuscule representations, without any
reference to Lie groups or Lie algebras.

Next we present an algorithm that constructs the
{\it fundamental adjoint varieties}   using the ideals of the
tangential and secant varieties of 
certain minuscule varieties.  By an {\em adjoint
variety},  we mean the unique closed orbit in the projectivization
$\PP\fg$ of a simple complex Lie algebra
$\fg$. We say that an adjoint variety is {\it fundamental} if the adjoint
representation is fundamental.   In particular, we construct all complex
simple Lie algebras without any reference to Lie theory.

Complex simple 
Lie algebras were first classified by Cartan and Killing 100 years ago.
Their classification proof proceeds by reducing the question to a
combinatorial problem: the classification of irreducible root systems,
and then classifying   root systems.

We present a new proof of the
classification of minuscule varieties and complex simple
Lie algebras by showing our
algorithms produce all minuscule (resp. fundamental adjoint) varieties
without using the classification of root systems, although
we do use some properties of root systems. We also provide
a proof that the only non-fundamental
adjoint varieties  are the adjoint varieties
 of $A_m$ and $C_m$, and thus we obtain a
new proof of the classification of complex simple Lie algebras.

\smallskip

Our proof can be translated into a combinatorial argument:
 the construction
consists of two sets of rules for adding new nodes 
to marked Dynkin diagrams. As a combinatorial algorithm, it is less
efficient than the standard proof, which proceeds by ruling out all but a
short list of Dynkin diagrams immediately, and then studying the few
remaining diagrams to see which are actually admissible.

\smallskip

Our constructions have applications that go well beyond 
the classification proof presented in this article.
This is the second paper in a series. In \cite{LMmagic}, \cite{LMseries} and
\cite{LM3} we   present geometric and representation-theoretic
applications of our algorithms.
In \cite{LMmagic} we   study the geometry 
of the exceptional homogeneous spaces using the constructions of this
paper. In \cite{LMseries} and \cite{LM3} we apply the results of this paper, 
especially our observations about the Casimir in section 5, to obtain
decomposition and dimension formulas for tensor powers of some
preferred representations. 

\section{Statements of main results}

Let $V$ be a complex vector space and let
 $X\subseteq\ppp V$ be a variety in the associated projective space.
Let $v_d(X)\subset\ppp (S^dV)$ denote its {\it $d$-th
Veronese re-embedding}. If $P_1\hd P_N$ is a basis of $S^dV^*$,
the space of homogeneous polynomials of degree $d$ on $V$, then the
map $\ppp V\ra\ppp(S^dV)$ is $[x]\mapsto [P_1(x) \hd P_N(x)]$. 

If $X\subseteq\ppp V$ and $Y\subseteq\ppp W$, we let
$Seg(X\times Y)\subset\ppp (V\ot W)$ denote their {\it Segre product},
given by $([x],[y])\mapsto [x\ot y]$. The Segre product generalizes to
an arbitrary number of factors.

We will use the notation
$\langle X\rangle \subset\ppp V$ to denote the linear span of $X$.

\begin{defi}\label{d1} Call a variety $X\subset\ppp V$ {\em a minuscule
variety} if $X=G/P_\a$ where
$G$ is a complex simple Lie group,
 $\a$ is a minuscule root, $P_{\a}$ is an associated maximal parabolic
subgroup and  $X$ is the projectivized orbit
of a highest weight vector in $V=V_{\o}$ where $\o$ is the fundamental weight dual
to the coroot of $\a$ (so the embedding is the minimal  equivariant embedding).
 Call $X$ a {\em generalized
minuscule variety} if $X$ is a Segre product of (Veronese re-embeddings of)
minuscule varieties. In this situation we will call $V$ a {\em minuscule} (resp. {\em generalized
minuscule}) $G$-module.\end{defi} 

The generalized minuscule varieties are those varieties admitting a
Hermitian symmetric metric induced from a Fubini-Study metric on the
ambient projective space. The minuscule varieties are those for which the
metric is irreducible and the embedding is minimal (i.e.,  not a Veronese
re-embedding).

\begin{defi}\label{d2}
 For a smooth variety $X\subset\ppp V$, let ${\mathcal
T}(X)\subset G(2,V)\subset \ppp
(\La 2 V)$ denote the {\rm   variety of  embedded 
tangent lines of $X$}. 

Let $\tau (X)\subset\ppp V$ denote the {\rm tangential variety of $X$} (the
union
of the points on embedded tangent lines), and let $\sigma_p(X)\subset\ppp
V$ denote the {\rm variety of secant $\pp{p-1}$'s to $X$}, that is, for
$\ux 1\hd \ux p\subset \ppp V$, let $\ppp_{\ux 1\hd\ux p}$ denote the
projective space they span (generally a $\pp{p-1}$), then
$\sigma_p(X)=\overline{\cup_{\ux 1\hd\ux p\in X}\ppp_{\ux 1\hd\ux p} }$. 
We let $\sx= \sigma_2(X)$.\end{defi}

\subsection{Minuscule case}

\noindent{\bf The minuscule algorithm}.
Let $Y=Seg( v_{d_1}(X_1)\times 
\cdots\times v_{d_r}(X_r))\subset\pp{n-1} =\ppp T$
where the $X_j\subseteq\pp {N_j}$'s are outputs of previous runs through
the algorithm or $\pp 1\subseteq\pp 1$.

We will call $Y$
{\it admissible} if ${\mathcal T}(Y)$
is linearly nondegenerate, that is, if
$\langle {\cT}(Y)\rangle =\ppp (\La 2 T)$. If $Y$ is admissible, then
define a rational map as follows: let  $d$ be the 
positive integer such that $\sigma_{d-1}(Y)\neq
\sigma_d(Y)=\pp{n-1}$. Linearly embed $\pp{n-1}\subset\pp n $ as
the hyperplane $\{x_0=0\}$,  and consider the rational map
$$\begin{array}{rccl}
\phi : & \pp n & \dashrightarrow & \pp N\subset\ppp (S^d\bcc{n+1*}) \\ 
& [x_0\hd x_n] & \mapsto &
[x_0^d,x_0^{d-1}T^*,x_0^{d-2}I_2(Y),x_0^{d-3}I_3(\sigma_2 (Y)) \hd
I_d(\sigma_{d-1} (Y))],
\end{array}$$
and call $X=\phi (\pp n)\subset \pp N$ an {\it output}.
 Here $T^*$ and $I_k(Z)=I_k(Z,\ppp T)$ 
 are shorthand notation  respectively for
a  basis of $T^*$ and a set of generators of the ideal of $Z$ in degree
$k$.

\begin{rema} We show below that the admissiblity hypothesis implies
$r,d_j\leq 2$
and thus there is a finite number of varieties to test at each run through
the algorithm. (Moreover, if $d_1=2$ then $r=1$.)\end{rema} 

\begin{theo}\label{t1} ({\rm Geometric construction of minuscule 
varieties}.) The
minuscule varieties are exactly the outputs of the minuscule
algorithm.\end{theo}

\begin{theo}\label{tt1} ({\rm The minuscule algorithm is effective}.) After six
runs the minuscule algorithm stabilizes and one can determine all minuscule
varieties from the output of the first six runs. In particular, the
minuscule algorithm gives an effective
classification of minuscule varieties.\end{theo} 

\subsection{Adjoint case}

\noindent{\bf The adjoint algorithm}.
Let $Y\subset\pp{n-2}=\ppp T_1$ be a
generalized minuscule variety. Here we define $Y$ to be
 {\it admissible} if $\langle {\cT}(Y)\rangle\subseteq\ppp (\La 2 T_1)$
has codimension one. Note that admissibility
implies there is an up to scale two-form defined on $T_1$. If $Y$ is
admissible, define a rational map as follows: linearly embed
$\pp{n-2}\subset\pp{n-1}\subset\pp n $ respectively as the hyperplanes
$\{x_n=0\}$ and $\{x_0=0\}$, consider the rational map
$$\begin{array}{rccl}
\phi : & \pp n & \dasharrow & \pp N\subset\ppp (S^4\bcc{n+1*}) \\ 
& [x_0\hd x_n] & \mapsto &
[x_0^4, x_0^3T_1^*,x_0^3x_n, x_0^2I_2(Y),x_0^2x_nT_1^*-x_0 I_3(\tau
(Y)_{sing}),x_0^2x_n^2-I_4(\tau (Y))], \end{array}$$
and call $X=\phi (\pp n)\subset \pp N$ an {\it output}. \smallskip

Note that to make sense of our notation in the mapping, we are using that
$\tau (Y)$ is a quartic hypersurface, which is proved in \S 6.1. In particular $I_3(\tau (Y)_{sing})$ is the space of derivatives of the equation
of $\tau (Y)$.

\smallskip

\begin{rema}We show below that the admissibility hypothesis implies that
$d_j, r\leq 3$, so again, there are a finite number of cases to check and
the algorithm is effective. (Moreoever $d_1=3$ implies $r=1$ and $d_1=2$
implies
$r\leq 2$.)\end{rema}

\begin{theo} 
The fundamental adjoint varieties are exactly the varieties constructed by
the adjoint algorithm. \end{theo}

\begin{coro} The fundamental adjoint varieties
are not rigid to order two. More precisely, there exist varieties
$X^n\subset\pp N$
  with the same projective second fundamental form as a fundamental adjoint
variety  at a general point which are not
fundamental adjoint varieties.\end{coro}

\begin{proof} If we write the image of $\phi$ as a graph in local
coordinates about $[1, 0\hd 0]$, the Taylor series will truncate
after the fourth order term. If we take the same Taylor series and
truncate after the second order term, we obtain, in this affine
open subset, a variety with the same second fundamental form  
as   $\phi (\pp n)$ at  $[1, 0\hd 0]$ on the entire
affine open subset, and this open variety
can be completed to a projective variety. (In general, varieties written
as graphs with only second order Taylor series have constant second
fundamental forms in region where the coordinates are valid.) On the other
hand, we show below that the image of $\phi$ is homogeneous
and thus has the same second fundamental form at all points.\end{proof}

\begin{rema} The corollary contrasts   the
case of the Severi varieties
and other  minuscule varieties of small codimension, which are
rigid to order two; see \cite{Lrigid}.\end{rema}

Our classification proof will be complete once we prove 
that the only simple Lie
algebras where the adjoint representation is not fundamental are $A_m$ and
$C_m$. We prove this in \S 6.3.

\begin{rema}
In many ways, the exceptional groups are better behaved than the classical
groups from the perspective of our constructions. For example:

-- The only \lq\lq exception\rq\rq\ in the minuscule algorithm is the
Grassmanian $X=G(k,V)$, constructed from $Y=\pp{k-1}\times\pp{l-1}$. In all
other cases the $P$-module $\La 2 T$ contains no $P$-submodules 
and there is no need to
study
$\langle \cT (Y)\rangle $.

-- The only pathological (i.e., non-fundamental) adjoint varieties are those
of $A_m, C_m$.

-- For the minuscule algorithm, the admissible $Y$'s which yield exceptional
$X$'s are such that $\tcodim\sigma (Y)$ is at most one, so the rational map
$\phi$ is at most cubic. \end{rema}

\section{Examples}
\subsection{Minuscule case}

For the first run through the algorithm, the admissible varieties and their
outputs are
$$\begin{array}{rclcrcl}
Y & \subseteq & \pp{n-1}& & X^n & \subseteq & \pp N\\ \pp 1 & \subseteq &
\pp 1 & &\pp 2 & \subseteq & \pp 2\\ \pp 1\times\pp 1 & \subset & \pp 3 &
&\QQ^4 & \subset & \pp 5\\ v_2(\pp 1) & \subset & \pp 2 & & \QQ^3 & \subset
& \pp 4. \end{array} $$
Here and below,
$\QQ^m\subset\pp{m+1}$ denotes the smooth quadric hypersurface. 

For the second round,
$$\begin{array}{rclcrcl}
Y & \subseteq & \pp{n-1} & & X^n & \subseteq & \pp N\\ \pp 2 & \subseteq &
\pp 2 & &\pp 3 & \subseteq & \pp 3\\
v_2(\pp 2) & \subset & \pp 5 & & G_{Lag}(3,6) & \subset & \pp {11}\\ \pp
1\times\pp 2 & \subset & \pp 5 & &G(2,5) & \subset & \pp 9\\ \pp 2\times\pp
2 & \subset & \pp 8 & & G(3,6) & \subset & \pp{19}\\ \QQ^4 & \subset & \pp
5 & & \QQ^6 & \subset & \pp 8\\
\QQ^3 & \subset & \pp 4 & & \QQ^5 & \subset & \pp 7 .\end{array} $$
Here $G(k,l)$ denotes the Grassmanian of $k$-planes in $\BC^l$ and
$G_{Lag}(k,2k)$ denotes the
Grassmanian of Lagrangian $k$-planes for a given symplectic form. 

Continuing, one gets
$$\begin{array}{rclcrcl}
Y & \subseteq & \pp{n-1} & & X^n & \subseteq & \pp N\\ \pp {n-1} &
\subseteq & \pp {n-1} & &\pp n & \subseteq & \pp n\\ v_2(\pp {m-1}) &
\subset & \pp {\binom {m+1}2-1 } & & G_{Lag}(m, 2m) & \subset & \pp
{C_{m+1}-1}\\ \pp {k-1}\times\pp {l-1} & \subset & \pp{kl-1} & & G(k,k+l) &
\subset & \pp {\binom{k+l}{k}-1} \\
\QQ^{2m-2} & \subset & \pp {2m-1} & & \QQ^{2m} & \subset & \pp{2m+1}\\
\QQ^{2m-1}& \subset & \pp {2m} & &\QQ^{2m+1} & \subset & \pp{2m+2}\\ G(2,m)
& \subset & \pp{\binom m2 -1} & & \BS_m & \subset & \pp{2^{m-1}-1}.
\end{array} $$

Here $C_{m+1}=\frac{1}{m+2}\binom{2m+2}{m+1}$ is the $(m+1)$-st Catalan
number, and the spinor variety $\BS_m$ of $D_m$ consists of one family of
maximal isotropic subspaces of $\CC^{2m}$ endowed with a nondegenerate quadratic form and embedded in the projectivization of one of
the two half-spin representations. 

Continuing, one sees that
$G_{Lag}(m,2m)$ and $\BS_m$ are terminal except for $\BS_5$ which yields
two exceptional spaces: $$\begin{array}{rclcrcl}
Y & \subseteq & \pp{n-1} & & X^n & \subseteq & \pp N\\ \BS_5 & \subset &
\pp{15} & & \BO\pp 2 & \subset & \pp{26}\\ \BO\pp 2 & \subset & \pp{26} & &
G_w(\BO^3,\BO^6) & \subset & \pp{55}.
\end{array}$$
Here $\BO\pp 2$ denotes the sixteen dimensional Cayley plane, which   is
the $16$-dimensional  homogeneous variety of
$E_6$, and $G_w(\BO^3,\BO^6)$ denotes the $27$-dimensional minuscule
variety of $E_7$.

\proof[Proof of \ref{tt1}]
First observe that $\langle \cT (Y)\rangle \neq \ppp (\La 2 T)$ for the
Segre product of three projective spaces. In fact, for $Y=Seg (\pp 1\times
\pp 1\times \pp 1)$, $\tcodim \langle\cT (Y)\rangle=1$ (and $Y=Seg (\pp
1\times \pp 1\times \pp 1)$ yields the $D_4$ adjoint variety), and $\tcodim
\langle\cT (Y)\rangle >1$ for all others. Thus one cannot have a triple
Segre product $Seg(Y_1\times Y_2\times Y_3)$ for any $Y_j$ in the minuscule
algorithm, and for any $Y_j$'s other than three $\pp 1$'s in the adjoint
algorithm.

For the Segre product of two varieties
$Seg (Y_1\times Y_2) \subset \ppp (W_1\otimes W_2)$, $\La 2 W_1\ot
I_2(Y_2)$ will not be in
$\langle\cT (Y)\rangle$ and similarly with roles reversed. (The ideals of
homogeneous varieties are generated in degree two; see \cite{lich}.) So
only a double
product of projective spaces is admissible for the first algorithm and in
the second algorithm, only $\pp 1\times \QQ^m$ is possible.

 For a Veronese embedding of a projective space, only $Y=v_2(\pp m)$ has
$\cT (Y)$ linearly full, and only $v_3(\pp 1)$ has $\langle\cT (Y)\rangle$
of codimension one, so no Veronese re-embedding of a subvariety of $\pp m$
will be admissible for either algorithm.

Thus, other than the examples mentioned above, the only new inputs to the
minuscule algorithm must be outputs from a previous round.
  From the list of examples above, we see that the algorithm stabilizes
after the sixth iteration. \qed

\subsection{Adjoint case}

The generalized minuscule varieties yielding fundamental adjoint varieties
are
$$\begin{array}{rclcc}
Y & \subset & \PP^{n-2}& &G\\
v_3(\pp 1)& \subset & \pp 3 & & G_2\\
\pp 1\times \QQ^{2m-3} & \subset & \pp {4m-3}& & B_m\\ \pp 1\times
\QQ^{2m-4}& \subset & \pp {4m-5}& & D_m\\ G_{Lag}(3,6) & \subset & \pp
{13}& & F_4\\ G(3,6)& \subset & \pp {19}& & E_6\\
\BS_6& \subset & \pp {31}& & E_7\\
G_w(\BO^3, \BO^6)& \subset & \pp {55}& & E_8 .\end{array}$$
The two exceptional (i.e., non-fundamental) cases are $$\begin{array}{rclcc}
\pp {k-1}\sqcup\pp {k-1} & \subset & \pp{2k-1} & & A_k\\ \emptyset &
\subset & \pp {2m-1} & & C_m . \end{array} $$

\begin{rema} From the examples one sees that the admissibility hypothesis
for the minuscule algorithm is equivalent to requiring that $\tau
(Y)=\sigma (Y)$ and admissibility for the adjoint algorithm is equivalent
to requiring that $\tau (Y) \neq \sigma (Y) =\ppp T$. Although it turns out
to be easier to use $\cT (Y)$, our constructions were motivated by two
theorems regarding
secant and tangential varieties: the Fulton-Hansen theorem \cite{FH}
 and Zak's theorem on Severi varieties \cite{Z}.
 The construction in the minuscule algorithm is a
generalization
of Zak's construction of the Severi varieties.\end{rema} 

\section{Local differential geometry}

In this section, we characterize the generalized minuscule varieties $G/P$
defined in \S 2 in terms of local differential geometry. We choose a Borel
subgroup $B$ of $G$ and a maximal torus $T$ of $B$, yielding a root system of
the Lie algebra
$\fg$, with a
base of the simple roots. We say that a simple root $\a$ is {\em short}
 if there exists a longer root inside the root system of $\fg$. Up to
conjugation, there is a natural correspondence between simple roots $\a$, 
 maximal parabolic subgroups $P_{\a}$ and nodes of the Dynkin diagram
of $G$.
We begin with general considerations.
\smallskip

Let $X^n\subset\ppp V$ be any variety and let $x\in X$ be a smooth point.
Let $\hat T_x\up k X\subseteq V$ denote the cone over the $k$-th osculating
space
and let $N_k=\hat T\up k/\hat T\up{k-1}$ denote the (twisted) $k$-th normal
space of $X$ at $x$ (see \cite{LM0}, \S 2.1). If we choose local
coordinates $(\xx 1\hd \xx N)$ at $x$ on $\PP V$, adapted to the filtration
by
osculating spaces, we may write $X$ locally as a graph. Taking
$\frac{\partial}{\partial\xx\a}$ as a basis of $T_xX$, for $1\leq \a\leq n$,
we have in these coordinates, for $\tdim\hat T\up {k-1}<\mu_k\leq \tdim\hat
T\up {k}$, the Taylor expansion:
$$\begin{array}{rcl}
\xx{\mu_2}&=& \qq{\mu_2}\a\b\xx\a\xx \b + \rr{\mu_2}\a\b\g\xx\a\xx \b\xx\g
+\hdots\\
\xx{\mu_3}&=& \qq{\mu_3}\a{\b\g}\xx\a\xx \b\xx\g + \rr{\mu_3}\a\b{\g\delta
}\xx\a\xx \b\xx\g\xx\delta +\hdots\\ &\vdots &
\end{array}$$
with summations over repeated indices.
The projective differential invariants of $X$ at $x$ are as follows: the
 {\em fundamental forms}
$\ff k =\fff k {X,x}=F^{k}_{k+0,X,x}\in S^kT^*\ot N_k$ 
are given by $$\ff 
k  =
\qq{\mu_k}{\a_1}{\hdots\a_k}d\xx{\a_1}\circ\hdots\circ d\xx{a_k} \ot
\frac{\partial}{\partial\xx{\mu_k}},$$ and the {\em relative differential
invariants} $F^{k}_{k+l, X}$ are given by 
$$F^{k}_{k+l,X,x }=\rr{\mu_k}{}{\a_1}{\hdots\a_{k+l}}d\xx{\a_1}\circ
\hdots\circ d\xx{\a_{k+l}}
\ot \frac{\partial} {\partial\xx{\mu_k}},$$ where the
$\rr{\mu_k}{}{\a_1}{\hdots\a_{k+l}}$ are the coefficients of the terms of
degre $k+l$ in the Taylor expansion of $x^{\mu_k}$ (see \cite{lanci} for
more precise definitions). Note that these
are a complete set of differential invariants in the sense that one can
recover the Taylor series of $X$ at $x$ from them. 

We use the notation $F_k=F_{k,X,x}^2$, $II=\fff 2 {X,x}$, $|II|=\fff 2 X
(N^*)\subseteq S^2T^*_xX$, and $\tbase |II| =\{ [v]\in \ppp T\mid
II(v,v)=0\}$.
\smallskip

For any variety, $F^{k}_{k+l,X}$ occurs as the coefficient of $\fff
{k+l}{v_{l+1}(X)}$ in the direction of $N_{k}\circ (\hat x )^{l}$ (with
respect to the induced adapted framing); see \cite{lanci}, 3.10. 

\medskip\noindent{\em Example.} Consider the cubic form $F_3=F_{3,X}^2$.
$F_3(Uv, Yv, Zv)$ is given by the component of $UYZv$ in $U(\fg)_2v \cap
N_2 $, that is, the component of $UYZv $ in $U(\fg)_2v$. Thus $F_3$
consists of sums of terms $UYZv$
such that $UYZv=\Sigma ^jW^j_1W^j_2v$ for some $W^j_1,W^j_2\in\fg$. For
example, if there are relations
in $U(\fg)$ that enable one to write $UYZ = W_1W_2$ for some
$U,Y,Z,W_1,W_2\in \fg$ such that $W_1W_2v\neq 0$, then $F_3\neq 0$.
\medskip

\begin{theo}\label{p1} Let $X\subset\ppp V$ be a homogeneous variety which
is a product of irreducible varieties of type $G/P_{\a}$, with $G$ simple
and $\a$ not short. Let $x\in X$. Then $X$ is generalized minuscule if and
only if the only nonzero differential
invariants of $X$ at $x$ are its fundamental forms.\end{theo} 

\medskip
\begin{rema} The result is false
if one drops the hypothesis of $X$ being homogeneous, even if one requires
$X$ to be smooth and $x$ to be a general point. However we believe it is
true if one does not require $X$ to be homogenous, but requires that the
only nonzero differential invariants are the fundamental forms at all
points of $X$. \end{rema}

\begin{rema}
If one thinks of the generalized minuscule varieties as those admitting a Hermitian symmetric
metric, then one could prove the result by observing that $\nabla
II^{herm}=0$, where $II^{herm}$ is the Hermitian second fundamental form
and $\nabla$ is the covariant differential operator. $II$, the projective
second fundamental
form, is the holomorphic part of
$II^{herm}$ and
$F_3$ is the holomorphic part of $\nabla II^{herm}$, and so on. Since we will
stay in the projective category, we will argue along different lines which
will give more precise information about the invariants $F_k$ in the other
cases. \end{rema}

\proof By \cite{lanci}, the differential invariants of Segre products and
Veronese re-embeddings can be computed 
from the original embeddings and the theorem holds in general  if it
holds for irreducible embeddings. Thus we restrict ourselves to the case
where $X=G/P_{\a}\subset\ppp V$ is irreducible and in its minimal 
embedding. 

Let $v\in V$ be a highest weight vector. Then   
  $T_{[v]}X\simeq \fg_+v$, where  
$$\fg_+ := \bigoplus_{\{ \b\in \Delta\mid (\o_{\a},\b)>0\}
  }\fg_{-\b}.$$
Here $\Delta$ denotes the set of all roots, $\o_{\a}$
denotes the weight dual to the coroot of  $\a$ and $(,)$
the invariant bilinear form.
Consider $X_1\cdots X_k v$, where $X_j\in\fg_+$. In order to 
have some nonzero component in $N_p$ with $p<k$,   
  some bracket among the $X_j$ must be nonzero. However, if $V$ is 
minuscule, then $[X,Y]=0$ for all $X,Y\in\fg_+$. This 
shows that $F^k_{k+l,X}=0$ for all $l>0$ when $X$ is minuscule. 

Conversely, let $Y=\tbase |II_{X,[v]}|$.  
 By, e.g., \cite{LM0}, $T_{[v]}X$ splits as a sum $T=T_1\op T_2\op
T_3\op\cdots$ of irreducible $H$-modules, where $H$ is a
maximal  semisimple
subgroup of $P_{\a}$, and $T=T_1$ if and only if $\a$
is minuscule.   

Let $Zv\in T_1$ be a general point of the cone over $\sigma (Y)$. We may
write
$Z=Z_1+Z_2$ with $Z_j^2v=0$, i.e., $Z_jv\in\hat Y$, the cone over $Y$.
We calculate
$$\begin{array}{rcl}
Z^2v & = & (Z_1Z_2+Z_2Z_1)v=(2Z_1Z_2 - [Z_1,Z_2])v, \\
Z^3v & = & (2[Z_1,Z_2](Z_1-Z_2) + (Z_1-Z_2)[Z_1,Z_2])v. 
\end{array}$$
 Note that if $Z^3v$ is not zero, it is in $N_2$ and thus 
$F_3(Zv,Zv,Zv)\ne 0$. We can write $Z^3v\equiv
3[Z_1,Z_2](Z_1-Z_2)v \;\;\tmod T_3$; thus 
we are reduced to showing that for
generic $Z$, $W_2W_1v\neq 0$ where $W_2=[Z_1,Z_2]\in T_2$ and
$W_1=(Z_1-Z_2)\in T_1$. 

If $X$ is not minuscule, i.e., if $\fg_2\simeq T_2\neq 0$, 
then $\fg_2= [\fg_1,\fg_1]$. In this case $[Z_1,Z_2]\neq 0$
because
  $Y\subset\ppp T_1$ is linearly 
nondegenerate and $Z_1,Z_2$ are general points of $\hat Y$.
  If $\a$ is not short, then 
$T^*_1\ot T^*_2\subset |II|$ by \cite{LM0}, 2.19, and thus
$W_2W_1v\neq 0$ for all $W_1\in T_1\backslash 0,\ W_2\in
T_2\backslash 0$. This proves our claim. \qed

\section{Proof of the minuscule case}

\noindent{\bf Strategy.}
We need to show that the outputs of the minuscule algorithm are
indeed minuscule varieties, and that all minuscule varieties
arise by the algorithm.

 In our algorithm, the rational map $\phi$ is defined such that at the
point
$x=[1,0\hd 0]\in X$, the only nonzero differential invariants
are the fundamental forms. Moreover,
$|\fff k{X,x}|= I_k(\s_{k-1}(Y))$  and $Y$ is generalized
minuscule. (Here and below, for an algebraic set $Z\subset\ppp V$,
 $I_k(Z)\subset S^kV^*$ denotes the component of the ideal of $Z$ in
degree $k$.)
By \cite{LM0}, 2.19 and 3.8, and   Theorem 4.1, 
any potential minuscule variety must be constructed out
of a generalized minuscule variety $Y$ by a mapping
of the form $\phi$.

We need to show that the additional hypothesis that
$\cT (Y)$ is linearly nondegenerate is necessary and sufficient
to imply that $X=G/P_{\a}\subset\ppp V$, where $G$ is simple and $\a$
is minuscule.

We proceed by constructing
the Lie algebra $\fg$, showing that
there is a unique candidate for $\fg$ and that this candidate
 can be given an appropriate Lie algebra structure
if and only if $\cT (Y)$ is linearly nondegenerate.

When such a $\fg$ exists, we then observe that the associated
minuscule variety has the same 
projective differential invariants as $X$ at a point and
therefore the two must coincide.
\smallskip

\noindent{\bf Analysis of a minuscule grading.}
A minuscule root  of a simple Lie algebra $\fg$  induces
  a  three step $\ZZ$-grading of $\fg$
$$\fg = \fg_{-1}\oplus\fg_0\oplus\fg_1,$$
 where $\fg_0=\fh\op\BC$ is a reductive Lie algebra with one
 dimensional center and having semisimple part $\fh$.
In addition $\fg_1$ 
can be identified with 
   $T= T_{[v]}G/P_{\a}$, where 
 $v\in V$ is  a highest weight vector. 
The closed orbit $Y=H/Q\subset\ppp T $ is a
generalized minuscule variety (see \cite{OV} and \cite{LM0}).

This $\ZZ$-grading of $\fg$ induces $\fg_0$-module
structures on $\fg_{\pm 1}$
which determine the brackets $[\fg_0, \fg_{\pm 1}]$. Moreover, $[\fg_1,
\fg_1]= [\fg_{-1}, \fg_{-1}]=0$. Thus   given  $\fg_{\pm 1}$ with 
their $\fg_0$-module structures, 
the only bracket we have not yet determined
is  $[\fg_1,\fg_{-1}]$. This bracket must be   a $\fg_0$-equivariant map 
$$\fg_{-1}\otimes \fg_1\longrightarrow \fg_0=\fh\op\CC.$$
This map has two components, the second of which induces a natural
duality between $\fg_1\simeq T$ and $\fg_{-1}\simeq T^*$, hence an
identification  of $\fg_{-1}$ with $T^*$. We denote the first component by
$\th$.

Consider the Killing form $B_{\fg}$ of $\fg$. Up to some nonzero 
constants, we must have, for $X_0,Y_0\in\fg_0$ and $X_{\pm 1}
\in \fg_{\pm 1}$, 
$$B_{\fg}(X_0,Y_0)=B_{\fh}(X_0,Y_0), \quad 
B_{\fg}(X_1,X_{-1})=\langle X_1,X_{-1}\rangle,$$
where this last bracket is given by the natural pairing between 
$\fg_1$ and $\fg_{-1}$. The invariance of the Killing form 
then implies that 
$$B_{\fh}(X_0,[X_1,X_{-1}])=B_{\fg}(X_0,[X_1,X_{-1}])=
B_{\fg}([X_0,X_1],X_{-1})=\langle [X_0,X_1],X_{-1}\rangle.$$
In particular, the map $\th$ is determined, up to a constant, by 
the $\fg_0$ action on $\fg_1$ and its pairing with $\fg_{-1}$. 
 The above identity can be rewritten, again up to a constant,
for $u^*\in \fg_{-1}$, $v\in \fg_1$, as 
$$\th (u^*\ot v) = \sum_i\langle u^*, X_iv\rangle Y_i,$$
where $X_i$ and $Y_i$ are dual bases of $\fh$ with respect to its
Killing form. \medskip

\noindent{\bf Construction of the Lie algebra.}
Now we are given a generalized minuscule variety $Y\subset\PP T$
for a semisimple Lie algebra $\fh$, and we want to  construct a
new Lie algebra $\fg$ and a   minuscule variety of $\fg$. In 
particular, $\fg$ must be endowed with a $\ZZ$-grading as above, and
we have just seen that the
only possible compatible Lie bracket in $\fg$ is   
determined, up to the multiplication of $\th$ by some constant, 
by the $\fh$-module structure of $T$. It remains to see whether
the constant can be chosen such that the Jacobi identities hold. 
\smallskip

Note that, assuming $\fg$ is a Lie algebra, it must be simple.
Indeed,   any nontrivial ideal $\fii\subseteq \fg$ 
is an $\fh$-module.
If the center $\BC\subset \fg_0$ were not in $\fii$, then
neither would   $T$ nor $T^*$ as $\BC\subset [T,T^*]$.
But since $[\fh, T]= T,[\fh, T^*]= T^*$, $\fii$ must be
zero in this case.  Now if $\BC\subset\fii$ we see
$T,T^*\subset\fii$ as well and thus $\fh\subset\fii$ so $\fii=\fg$.
\medskip
 
\begin{prop} \label{great}
The bracket defined above  endows $\fg$ with
the structure of a  Lie algebra for some constant multiple of
$\th$   if and only if 
${\cT}(Y)\subset\PP(\we^2T )$ is linearly nondegenerate.
\end{prop}

The choice of constant for $\th$ is unique; see the proof
of Lemma 5.2 below.

\begin{proof}
To determine if $\fg$ is actually a Lie algebra, we need
to check the Jacobi identity. This identity can be split
into a number of graded parts, which are of
different natures. The $(\fg_0,\fg_0,\fg_i)$ identities follow because
$\fg_i$ is a $\fg_0$-module, while the $(\fg_0,\fg_i,\fg_j)$ identities
follow because the bracket $[\fg_i,\fg_j]\rightarrow \fg_{i+j}$ is 
$\fg_0$-equivariant. It remains to verify the $(\fg_i,\fg_j,\fg_k)$ 
identities for $i,j,k\ne 0$. By symmetry, the only
identity     to check is the one involving $v,w\in \fg_1=T$ and 
$u^*\in\fg_{-1}=T^*$. We must show that
$$[[u^*,v],w] + [[v,w],u^*] + [[w,u^*],v] = 
\th (u^*\ot v)w -\th (u^*\ot w)v 
+\langle u^*,v\rangle  w-\langle u^*,w\rangle  v = 0,$$
i.e., that
$$t^* (\th (u^*\ot v)w -\th (u^*\ot w)v )= 
\langle t^*,v \rangle \langle u^*,w\rangle 
-\langle t^*,w \rangle \langle u^*,v\rangle 
$$
for all $t^*\in T^*$.
The map $(v,w,u^*,t^*)\mapsto  
\langle t^*,v \rangle \langle u^*,w\rangle 
-\langle t^*,w \rangle \langle u^*,v\rangle$,
when considered as an element of $\La 2 T^*\ot \La 2 T$, 
is just twice the identity, so we must prove that  
$\TH \in \La 2 T^*\ot \La 2 T$ defined by 
$$
\TH (v\ww w)(u^*,t^* )= t^* (\th (u^*\ot v).w -\th (u^*\ot w).v )
$$
is a homothety, since $\th$ is only determined up to a constant
multiple. 

\begin{lemm}  
   $\Theta \mid_{\langle\cT (Y)\rangle}$ is a homothety.
\end{lemm}

\begin{proof} Given a Cartan subalgebra $\ft$ of $\fh$ and 
a basis of the corresponding root system $\Delta$, 
we choose an adapted basis $E_\b, H_j, E_{-\b}$ of $\fh$, where $\b\in
\Delta_+$ is a positive root, $B_{\fh}(E_{\b},E_{-\b})=1$, 
and the $H_j$ give an orthonormal basis 
of $\ft$ with respect to the Killing form.  Then 
$$\TH (v\ww w) = 
\sum_{\b\in\Delta_+} (E_{\b} v\ww E_{-\b}w
+E_{-\b} v\ww E_{\b}w ) 
+ 2\sum_j  H_j v\ww  H_jw .$$

Now assume, to simplify notation, that $Y$ is a minuscule  variety 
 of the  simple group $H$, so that 
$Y=H/Q_{\zeta}$ for some parabolic subgroup $Q_{\zeta}$ of $H$ 
corresponding to a   minuscule root $\zeta$ of $\fh$. Then $T$ is the 
fundamental $\fh$-module of highest weight the fundamental weight 
$\o$ dual to the coroot of $\zeta$.  If $Y$ is a Segre product of
minuscule varieties, the argument will be unchanged. If $Y$ is
a Segre product of   Veronese re-embeddings
of minuscule varieties, it
is necessary to keep track of the degrees of the re-embeddings in the
argument. 

In the computation above, we may suppose that $v\in Y$ and $w$ is 
tangent to a
  line of $Y$ passing through $v$
  (as $\cT (Y)$ being linearly nondegenerate
  implies $Y$ is linearly nondegenerate). Since $Y$ is $H$-homogeneous, we may
assume that
$v$ is a highest weight vector for $T$. In this case, when $\b$ is a
positive root, $E_{\b}v=0$,  and $E_{-\b}v=0$  unless $\b$ has coefficient 
one on the minuscule root $\zeta$. By linearity, we may 
 assume that $w=E_{-\g}v$ for some positive root $\g$. Then, 
for $\b$ as above, $E_{\b}w=0$ unless $\b=\g$, and thus
$$\begin{array}{rcl}
\TH (v\ww wv ) &=&
  E_{-\g} v\ww E_{\g}w
 + 2\Sigma_j H_jv\ww H_jw
\\ &=& 
w\ww E_\g E_{-\g}v +2\Sigma_j\o(H_j)(\o-\g)(H_j)v\ww w
\\ &=& 
w\ww [E_\g,E_{-\g}]v +2(\o, \o-\g)v\ww w.
\end{array}$$
Here we denoted by $(\;,\;)$ the pairing on $\ft^*$ dual to the 
Killing form.
We have
$[E_\g, E_{-\g}]=\frac{(\g,\g)}{2}H_\g$, where $H_{\g}$ is 
the coroot of $\g$ (see   \cite{bou}). Finally,
$$\o ([E_{\g},E_{-\g}])=\frac{(\g,\g)}{2}\o
(H_{\g})=(\o,\g)=(\o,\zeta)$$
because the Killing form identifies
the  coroot of $\g$ with $\frac{2\g}{(\g,\g)}$ and  
the last equality holds because $\g$ has coefficient one on the simple 
root $\zeta$. The lemma is proved.  \end{proof} 

\begin{rema}It turns out that to accommodate semisimple Lie algebras, it will
be necessary to normalize long roots to have 
the same length independent of the Lie algebra.
We choose length two, so one obtains
$\TH|_{\langle\cT(Y)\rangle} = 2((\o,\o)-2)Id$.\end{rema}

\medskip For the opposite direction, we introduce a definition. 

\begin{defi} Let $\fh$ be a semisimple Lie  algebra and
let $U$ be an  $\fh$-module.
We will say $U$ is $C$-{\em irreducible} if the Casimir
operator $C_U$ acts on $U$ as a homothety. In particular, every
irreducible $\fh$-module is $C$-irreducible. 
If $U$ is $C$-irreducible, we let  $c_U$ 
be the constant such that $C_U=c_UId_U$. 

We say that $U$ is {\em almost} $C$-{\em irreducible} if $U=\CC\op V$
where $V$ is $C$-irreducible.
\end{defi}\smallskip

\begin{lemm}Let $\fh$ be   semisimple,
let $T$ be an irreducible
$\fh$-module and let $\Theta$
be defined as above.
Then $\Theta=C_{\we^2T }-2c_{T }Id_{\we^2T }$. In
particular,
$\TH$ is a homothety if and only if $\La 2 T $ is $C$-irreducible.
\end{lemm}

\proof
Let $X_i, Y_i$ be bases of $\fh$ dual with respect to the Killing
form. Let $v,w\in T $. Then
$$
\begin{array}{rcl}
C_{\we^2T }(v\ww w)&=&\Sigma_i X_iY_i(v\ww w)\\
&=&\Sigma_i (X_iY_iv)\ww w +  (Y_iv)\ww X_i w + (X_iv)\ww Y_iw
+ v\ww (X_iY_iw)\\
&=&2c_{T }v\ww w + \TH (v\ww w),
\end{array}$$
which implies our claim. \qed

\medskip
Our proof of Proposition \ref{great} will be complete once we
have proven the following lemma:
 
\begin{lemm} Let $T$ be a generalized minuscule module  of a 
semisimple Lie algebra 
$\fh$. Let $Y\subset\ppp T$ be the 
corresponding generalized minuscule variety.
If   $\we^2T$ is $C$-irreducible, 
then ${\cT}(Y)$ is linearly nondegenerate.
\end{lemm}

We first prove the case where $\fh$ is simple,
which follows from the following lemma:

\begin{lemm}\label{simple} Let
$\fh$ be a simple Lie algebra and let $T$ be a minuscule $\fh$-module  or 
  a symmetric power of a minuscule $\fh$-module. 
  Suppose that $\we^2T$ is $C$-irreducible. 
Then  $\we^2T$ is irreducible. \end{lemm}

\begin{proof}
Let $\o=\o_{\a}$ be the highest weight of the minuscule module $T$.  
The highest weight of $\we^2T$ is $2\o-\a$. 
Suppose that $\we^2T$ has an irreducible component of another
highest weight. This weight $\m$ must be the sum of two distinct
weights of $T$. Hence, since $T$ is minuscule, $\mu =u\o+v\o$ 
for two distinct elements $u,v$ of the Weyl group. We show
that $c_{\mu}<c_{2\o-\a}$, which will prove that $\we^2T$ is 
not $C$-irreducible. Since $(\a,2\r-\a)=0$, we have 
$$
c_{2\o-\a}-c_{\mu}=2(\o,\o)-2(u\o,v\o)+(2\o-u\o-v\o,2\r).
$$
Write $u\o=\o-\sum_{\g}n_{\g}(u)\g$, where the sum is over 
the simple roots. The coefficients $n_{\g}(u)$ are non-negative, and 
$n_{\a}(u)>0$ when $u\o\ne \o$. Thus 
$$\begin{array}{rcl}
(u\o,v\o) & = & (\o,u^{-1}v\o)=(\o,\o-\sum_{\g}n_{\g}(u^{-1}v)\g)
=(\o,\o)-n_{\a}(u^{-1}v), \\
(u\o,2\r) & = & (\o-\sum_{\g}n_{\g}(u)\g,2\r) = 
(\o,2\r)-\sum_{\g}(\g,\g)n_{\g}(u).
\end{array}$$
This implies that 
$$c_{2\o-\a}-c_{\mu}=\sum_{\g}(\g,\g)(n_{\g}(u)+n_{\g}(v))
+2n_{\a}(u^{-1}v)-4.$$
This number is positive: this is clear if $u\o,v\o\ne\o$, 
since $n_{\a}(u), n_{\a}(v), n_{\a}(u^{-1}v)$ are then all
positive; otherwise, we may suppose that $u=id$, in which case
$$c_{2\o-\a}-c_{\l}=\sum_{\g\ne\a}(\g,\g)n_{\g}(v))
+4n_{\a}(v)-4,$$
with $n_{\a}(v)>0$, and this number is clearly positive unless
$v\o=\o-\a$, which is excluded. 

The case of symmetric powers of a minuscule module is  
similar. \end{proof}

\smallskip

We consider now the case where $\fh$ is not simple, and $T=T_1\otc T_r$,
where $T_i$ is a $\fh_i$-module of the form $S^{d_i}U_i$, with $U_i$
fundamental and minuscule. For $r=3$, we have 
$$\we^2T=\we^2T_1\ot S^2T_2\ot S^2T_3\op S^2T_1\ot \we^2T_2\ot S^2T_3
\op S^2T_1\ot S^2T_2\ot \we^2T_3\op \we^2T_1\ot\we^2T_2\ot\we^2T_3.$$
For such a module to be $C$-irreducible, we first need that each 
$S^2T_i$ and $\we^2T_i$ be $C$-irreducible. But we always have 
$c_{\we^2T_i}<c_{S^2T_i}$, hence a contradiction. For the same reasons
we cannot have $r\ge 4$. 

If $r=2$, then $\we^2T=\we^2T_1\ot S^2T_2\op S^2T_1\ot \we^2T_2$
is $C$-irreducible if and only if the $S^2T_i$ and $\we^2T_i$ are
$C$-irreducible. To see this, recall that for a   representation
$V_{\o}$, if $S^2V$ and
$\La 2V$
are $C$-irreducible, then
 with the normalization 
of long roots having length two,
 $c_{S^2V}= c_{\La 2V}+ 4$, so  
$c_{\La 2T_1\ot S^2T_2}=c_{S^2T_1}-4+c_{S^2T_2}=c_{S^2T_1\ot \La 2T_2}$.
By the above lemma, this implies that $\we^2T_i$
  and   $S^2T_i$ are irreducible. But then
$\cT (Y)$ is linearly nondegenerate. 
\end{proof} 

\medskip
We record the following lemma, which we will need in \S 6.

\begin{lemm} If $\we^2T$ is almost $C$-irreducible, then
 $\langle {\cT}(Y)\rangle$ has codimension one.\end{lemm}

\begin{proof}
The lemma follows from 
Lemma \ref{simple} when $\fh$ is simple. For $T=T_1\otc T_r$, 
we check as above that $r\ge 3$ is not possible unless $r=3$ and 
the $T_i$ are two-dimensional. Finally, for $r=2$, we can
suppose that $S^2T_1$ and $\we^2T_2$ are $C$-irreducible, while
$\we^2T_1=\CC\op A_1$ and $S^2T_2=\CC\op A_2$ are $C$-irreducible.
But then $\we^2T$ contains $A_1$ and $A_1\ot A_2$, and the Casimir 
operator must act on both modules as the same constant. This leads
to a contradiction unless $A_1$ is zero, which means that $T_1=\CC^2$ is 
two-dimensional, and $\we^2T=S^2T_2\op S^2\CC^2\ot \we^2T_2$.
We note that in this case  $\langle\cT(Y)\rangle $ must 
have codimension one.
\end{proof}

\begin{rema} The minuscule and adjoint algorithms can
be generalized to construct all $G/P$'s, however the resulting
rational maps get more and more complicated. The corresponding
construction of a $\ZZ$-graded Lie algebra $\fg$ from a Lie algebra
$\fh$ equipped
with a representation $T_1$
 (by gluing in a node to the marked Dynkin diagram) is easier to
 describe and even provides an effective algorithm for producing $\ZZ$-graded
 affine Lie algebras.

 One takes $\fg_0=\fh+\CC$ and $\fg_1=T_1$ as before. Next, for
 $\fg_2$, one takes the $\fh$-module
 complementary  to $\langle \cT(Y)\rangle$ in $\La 2T_1$, and denote
 the corresponding closed $H$-orbit $Y_2\subset\ppp\fg_2$.
 Now there are further Jacobi identities to check, but these can
 be checked geometrically by considering the linear span of
 the natural incidence variety $Y_{12}\subset Y\times Y_2\subset
 \ppp S_{21}(T_1)$. If it is linearly full, one has
 already constructed $\fg$. If not,
 one takes the  complement of its linear span as $\fg_3$. At each
 step one needs only to consider the linear span of
   the analogously defined
  $Y_{1k}$. (This is because N. Mok has observed that in a $\ZZ$-graded
  Lie algebra generated by $\fg_1$, one always has
  $[\fg_p,\fg_q]\subseteq [\fg_1,\fg_{p+q-1}]$ (personal communication).)
   By way of example, consider
 the $\ZZ$-gradings of $\fe_8\up 1$ ($\ZZ_m$-gradings of
 $\fe_8$) obtained  from the extremal
 nodes: from $\a_1$ we have a $\ZZ_2$-grading with $\fg_0=\fd_8+\CC$,
 $\fg_1=V_{\o_8}$ (spin representation), $\a_2$ yields a $\ZZ_3$-grading 
   with $\fg_0=\fa_8+\CC$, $\fg_1=V_{\o_3}$, $\fg_2=V_{\o_6}$, and
 $\a_8$ yields a $\ZZ_2$-grading with $\fg_0=\fe_7+\fa_1+\CC$,
 $\fg_1=V_{\o_7}^{\fe_7}\ot V_{\o_1}^{\fa_1}$.
 
 One could attempt to continue the algorithm for general
 Kac-Moody Lie algebras, but since there will be exponential growth
 of the factors, the algorithm is not at all effective. Note that
 for affine Lie algebras the algorithm is considerably more explicit than the
 standard method of construction by taking the free Lie algebra and 
 then quotienting out by the relations.
 \end{rema} 

\subsection{Minuscule algorithm: the algebraic version} 

We may rephrase the minuscule algorithm algebraically: Let 
$\fh = \fh_1\op\cdots\op
\fh_r$ where the $\fh_j$ are constructed in an earlier run
through the algorithm. Let $T_j$ be the corresponding
representations constructed or a symmetric
power of such representations, and let $T=T_1\ot\cdots\ot
T_r$. Let $\tilde Y \subset\ppp T^*$ be the closed $H$-orbit, let
$V_0=\BC$, $V_1=T$, 
$V_2=I_2(\tilde Y )\subseteq S^2T$,  and let 
$$V_k=(V_{k-1}\ot T) \cap S^{k}T,$$
the {\em prolongation} of $V_{k-1}$ (see \cite{LM0}). 

Call $(\fh, T)$ {\it admissible} if $\langle \cT (\tilde Y )\rangle
=\ppp\La 2 T^*$. If $(\fh, T)$ is admissible, then
$$\begin{array}{rcl}
\fg & := & T^*\oplus (\fh\op\BC )\oplus T, \\
V & := & \oplus_{j=0}^{\infty}V_j
\end{array}$$
are respectively a simple Lie algebra, and a minuscule $\fg$-module.
(The last sum is finite because the terms correspond to the
generators of the ideals of successive secant varieties.) 

Moreover, starting with $(\fh,T)=(0,\BC)$ and iterating the
algorithm, one arrives
at all pairs $(\fg, V)$, with $\fg$ simple and $V$ a
minuscule $\fg$-module.

\begin{rema} Instead of constructing $\fg$ first, one could first
construct $V$ as above, define an action of each component of $\fg$
on $V$, and define the bracket in $\fg$ via the commutator
of the actions on $\fg$. The algebra $\fh$ acts naturally on each $V_j$,
one
makes $\BC\subset\fg_0$ act on each factor by a scalar, and $T^*$
acts naturally as a lowering (or \lq\lq annihilation\rq\rq ) operator
$V_j\ra V_{j-1}$ by differentiation. One makes $T$ act as a
raising (or \lq\lq creation\rq\rq ) operator as follows:  for
$p\in V_{l-1}$, one defines $w .p$, up to a constant, 
as the projection of $w\circ p$ on $V_l$ (one must check that this 
is well defined).

For example, consider the case $Y=G(2,W)$, $X=\BS$. Here
$V_k = \La{2(k+2)}W$
and the raising and lowering action corresponds to Clifford
multiplication. Thus we see that minuscule
representations define
algebraic structures  that are cousins of Clifford algebras.

Taking $V_2$ to be the ideal of a variety ensures the
sum is finite. To study infinite dimensional representations,
take $V_2\subseteq S^2T$ to be an $\fh$-invariant subspace 
whose successive prolongations do not terminate, e.g., one
could take $V_2$ to be involutive (see \cite{bcg} for the definition
of involutivity).
\end{rema}
 

\section{Proofs of the adjoint theorems}

\subsection{The fundamental case}

\noindent {\bf Strategy}. We need to prove that
  the adjoint algorithm is well defined, i.e., that
$\tau (Y)$ is a quartic hypersurface,
that the adjoint algorithm    constructs all fundamental
adjoint varieties, and that
  the adjoint algorithm only constructs fundamental adjoint 
varieties.

\smallskip

We   begin, as in the minuscule algorithm, by studying
the five step $\ZZ$-grading of a simple Lie
algebra $\fg$ whose adjoint representation is fundamental.
By \cite{LM0,OV}, 
$\fg_0=\fh\op\BC$, $\fg_1=T_1$ is a minuscule $\fh$-module
and $\fg_2=\CC$.
Letting $Y_1\subset\ppp T_1$ denote the closed $H$-orbit,
we   show that a necessary and sufficient condition that the pair
$(\fh, T_1)$ produces a simple Lie algebra $\fg$
together with its adjoint representation is   that
the linear span of $\cT (Y)$ has codimension one. (Note that 
the line in $\La 2 T_1$,
generated by the symplectic form of $T_1$,  already implies that the
complement of $\langle\cT (Y)\rangle$ must contain a trivial
$\fh$-module.)  

We   then determine the differential invariants of a fundamental
adjoint variety, in particular, we show  that $\tau (Y)$ must
be a quartic hypersurface.
 We finally observe that the  invariants
agree with the differential
invariants of the varieties $X$ constructed by the algorithm
at our preferred point $[1,0\hd 0]\in X$. \medskip

\noindent {\bf Analysis of the adjoint grading}. 
Let $\fg$ be a simple Lie algebra whose adjoint
representation   is fundamental. We choose a Cartan 
subalgebra of $\fg$  and a basis of the corresponding root 
system $\Delta$. We denote the highest root by $\tilde\a$, 
and consider the corresponding five step grading (\cite{OV})
$$\fg= \fg_{-2}\oplus \fg_{-1}\oplus \fg_{0}\op \fg_{1}\op\fg_{2}.
$$
Here $\fg_0= \fh\op\BC$, $\fg_2=\BC$, and $\fg_{1}=T_1$ 
is a generalized minuscule $\fh$-module. $T_1$ is the contact
hyperplane inside the tangent space to
a point     $x\in G/P_{\tilde\a}\subset\ppp\fg$
such that $x=\ppp\fg_{-2}$.
Let $Y\subset\ppp T_1$ be the closed $H$-orbit, then
  $Y=\tbase |II_{G/P_{\tilde\a},x}|$. 
(See \cite{LM0}.)

The Lie algebra structure of $\fg$ can be recovered  
as follows:   the Killing form of $\fg$
induces an identification between $\fg_{-1}$ and $T_1^*$,  and
also between $\fg_{-2}$ and $\fg_2^*$. 
The bracket 
$[\fg_1,\fg_1]\rightarrow \fg_2=\CC$ defines a symplectic form on
$T_1$. We   choose $\o^*\in\fg_{-2}\subset\La 2T_1$ and 
$\o\in\fg_{2}\subset\La 2T_1$
dual with respect to the Killing form. By contraction, $\o $ induces 
an identification of $T_1$ with $T_1^*$, and 
we let $u^*=\o (u,.)$. Then 
$$\begin{array}{rcll}
\hspace{0mm} [u,v] & = & p\o^*(u,v)\o &  \\
\hspace{0mm} [u^*,v^*] & = & p^*\o^*(u,v)\o^* & \forall u,v\in T_1, \\
\hspace{0mm} [u^*,\o] & = & mu, & \\
\hspace{0mm} [u,\o^*] & = & m^*u^* & \forall u\in T_1,
\end{array}$$
for some constants $p,p^*,m,m^*\in\CC$. By the invariance
of the Killing form,  
$$p\o^*(u,v)=B_{\fg}([u,v],\o^*)=B_{\fg}(u,[v,\o^*])
=m^*\langle u,v^*\rangle = m^*\o^*(u,v),$$
hence $p=m^*$ and similarly, $p^*=m$. 

There exist  $l,o\in\CC$ and an element $\11$ in
the center of $\fg_0$ such that
$$\begin{array}{rcl}
[u^*,v] & = & \o^*(u,v)\11+l\th(u^*\ot v) \quad\forall 
u,v\in T, \\ \hspace{0mm} [\o^*,\o] & = & o\11 .
\end{array}$$
Note that $\11$ must act on $T_1$ by some constant $a$. 
The map $\th$ must be
symmetric in the sense that 
$$\th(u^*\ot v)= \th(v^*\ot u)\qquad \forall u,v\in T_1$$
because   $\fh$ preserves the form $\o$.
The Jacobi identities give additional constraints on the 
scaling  constants $l,m,m^*,o,a$, and one can check that  
  $o=a=m=-m^*=2$. Moreover,
we will see that $l$ is determined by the Jacobi identity.
Thus the Lie algebra structure of $\fg$ is      determined 
by the $\fh$-module structure of $T_1$. 
\medskip

\noindent {\bf Construction of the Lie algebra}.
Now we  begin  with an $\fh$-module 
$T_1$ and an $\fh$-invariant two-form $\o$ on $T_1$ (defined up
to scale) and attempt to define a
Lie  algebra structure  on
the vector space $\fg$ with its five step gradation as above. 
The bracket on $\fg$ is determined up to the constant $l$ 
and we must check if the Jacobi identities hold. 

\begin{prop}
Let $\fh$ be semisimple and let $T_1$ be a generalized minuscule
$\fh$-module. Then there exists a choice of $l$
such that  $\fg$, with the bracket 
  defined  above, 
is a   Lie algebra  if and only if
$\tcodim\langle \cT (Y)\rangle =1$. In this case, $\fg$ is simple and its
adjoint representation is fundamental. 
\end{prop}

\begin{proof} 
By the same argument as in the minuscule case, the 
only identity we need to verify is 
the $(\fg_1,\fg_1,\fg_{-1})$ identity, that is, for all $u,v,w\in T_1$, 
$$l[\th(w^*\ot v)u-\th(w^*\ot u)v]
=\o^*(v,w)u-\o^*(u,w)v-2\o^*(u,v)w.$$
As before,  we interpret  this
equation as an equality between endomorphisms 
of $\La 2T_1$. The left hand side represents the action of the Casimir 
operator (plus some multiple of the identity), the first two terms of
the right hand side represent twice the identity of $\La 2T_1$, and 
the last term is the projection onto $\CC\o\subset\La 2T_1$. Thus
  the identity above will hold, with an appropriate
choice of $l$  if and only if $\La 2T_1$ is almost $C$-irreducible. 
By Lemmas 5.2 and 5.7, this is true if and only if 
$\tcodim\langle \cT (Y)\rangle =1$. 

One sees that $\fg$ is simple by an   argument 
similar to the minuscule case.\end{proof}

\medskip\noindent {\bf The adjoint variety}. 
We now calculate the differential invariants of a fundamental 
adjoint variety.
We   consider the adjoint variety $G/P_{\tilde\a}\subset 
\ppp \fg$ as  the orbit of $[\o]$, and identify the tangent space of 
$G/P_{\tilde\a}$ at $[\o ]$ with $\CC\op T_1\subset\fg_0\op\fg_1$.  
Let $X=x\o^*+u^*$ denote a general element of $\fg_{-2}\op\fg_{-1}$, 
which implies $X\o=2x\11+2u$ is a general element of
$T_{[\o ]}G/P_{\tilde\a}$ and $G/P_{\tilde\a}$
is  the projectivization of the closure of the
union of elements of the form $\exp
(X)\o$. We compute 
$$\begin{array}{lcl}
X\o &=& 2x\11+2u,\\
X^2\o &=& 8x^2\o^*+2l\th (u^*\ot u), \\
X^3\o &=& -2l\th (u^*\ot u)u^*, \\
X^4\o &=& 2l\; \o^*(u,\th (u^*\ot u)u), \\
X^5\o &=& 0.
\end{array}$$

This computation already shows that our map $\phi$ in the adjoint
algorithm must indeed be of degree (at most) four. It remains to 
interpret  the above equations. 
  From the second line and the 
fact that by \cite{LM0}, $\tbase |II|=Y$, we obtain:

\begin{prop}   $Y=\tbase |II|\subset\PP T_1$ 
is as follows:
$$
\begin{array}{rcl}
Y&=&\ppp \{  u\in T_1\mid  \th (u^*\ot u)=0 \}\\
&=& \ppp \{  u\in T_1\mid  \o^*(u,Zu)=0 \ \forall Z\in\fh \}.
\end{array}
$$ 
In particular, the second description implies that
$u$ is 
$\o^*$-orthogonal to the tangent space of its $H$-orbit.
\end{prop}

\begin{rema}In general, the space of quadrics vanishing on a closed orbit
$Y$ is the kernel of the natural map from $S^2T_1^*$ to $S^{(2)}T_1^*$,
where
$S^{(2)}T_1^*$ denotes  the Cartan product of $T_1^*$ with itself. An
important  consequence  of  Proposition 6.2 is  that 
$$S^2T_1=S^{(2)}T_1\op\fh.$$\end{rema}

Finally, we prove that the terms of degree three and four   
above are the same as those in the adjoint algorithm:

\begin{prop}  The tangential variety   
$\tau (Y)\subset\ppp T_1$ is the quartic hypersurface
defined by the equation 
$$p(w)=\o^*(\th (w^*\ot w)w,w)
=B_{\fh}(\th (w^*\ot w),\th (w^*\ot w))=0.$$ 
In particular, the singular locus of $\tau (Y)$ is given by the 
equations $\th (w^*\ot w)w=0$.
\end{prop}

\begin{proof}
   $\tau(Y)$ is nondegenerate because $III_Y\ne 0$, 
(\cite{lanri}, 10.1), hence it is a  hypersurface. 

We check that $w=v+Xv\in \tau (Y)$ implies  $p(w)=0$. 
Introduce the notation  $\th^*(uv)=\th^*(u.v)=\th(u^*\ot
v)$ and note that $\th^*(uv)=\th^* (vu)$.  For all $v\in Y$, we have $\th^*
(v^2)=0$. Differentiating, we obtain 
$\th^*(v.Xv)=0$ and  $\th^*((Xv)^2)+\th^*(v.X^2v)=0$ for all $X\in\fh$. 
Thus $\th^*(w^2)=\th^*((Xv)^2)$.  

Now   recall the identity 
$$l[\th^*(rs)t-\th^*(rt)s]=\o^*(r,s)t-\o^*(r,t)s+2\o^*(s,t)r,$$
which was a consequence of the Jacobi identity on $\fg$ (see the 
proof of Proposition 6.1).
Setting $r=s=Xv$ and  $t=v$, we obtain $l\th^*((Xv)^2)v =l\th^*(v.Xv)v+
 \o^*(Xv,v)Xv=0$.
Similarly, setting $r=v$, $s=X^2v$, and $t=Xv$, we obtain
$l\th^*(v.X^2v)Xv = 2\o^*(X^2v,Xv)v$. This implies that 
$\th^*(w^2)w = -2\o^*(X^2v,Xv)v$, so that finally 
$$ p(w)=\o^*(\th^*(w^2)w,w) = -2\o^*(X^2v,Xv)\,\o^*(v,Xv) = 0.$$
$p$ must be irreducible as
there is no invariant linear or quadratic form
on $T$, so
  $p$ must be   the equation of $\tau(Y)$. \end{proof}

\begin{rema} It turns out that the smooth part of the 
tangential variety is a single $\fg$-orbit. We have no 
{\it a priori} proof of this observation.
Such a proof would make the above
computation unnecessary. A  consequence of the observation is that 
the invariant ring of $T$ is free, and generated by the
quartic polynomial $p$.
\end{rema}

\subsection{Another perspective on the adjoint algorithm} 

We may think of
$\fg_{-2}\op\BC\op\fg_2\subset \fg$ as a copy of $\fsl_2$, where $\BC$
is the center of $\fg_0$. Thus we may write $\fg = T_1\op\, 
(\fh\oplus\fsl_2)\,\op T^*_1$ and 
think of our modified
construction as a \lq \lq quantization\rq\rq\ (in the sense that 
$\BC \ra \fsl_2$) of the construction in the minuscule case.

If we let $\BC^2 $ be the standard $\fsl_2$ module, we can
write $\fg$ more concisely (recalling that $T_1$ is a symplectic module,
hence can be identified with its dual) as 
$\fg=(\fh\oplus\fsl_2)\;\op\,\BC^2\ot T_1$.
All our constructions 
can be made using invariant tensors on these two spaces. 
In this way we recover the $\ZZ_2$-gradings of a simple
Lie algebra obtained by marking the node(s) of
the Dynkin diagram of the affine Lie algebra $\fg\up 1$
adjacent to $\tilde\a$. Here $\fg_{\underline 0}= (\fh +\fsl_2)$,
and $\fg_{\underline 1}= \BC^2\ot T_1$. In
particular, Proposition 6.1 can be (partly) rephrased as follows:

\begin{prop} Let $\fh$ be a simple Lie algebra and $T$ an irreducible 
symplectic $\fh$-module, such that $\La 2T$ is the sum of an
irreducible module and a line. Then 
$$\fg=(\fh\oplus\fsl_2)\,\op\,\BC^2\ot T$$
has a natural structure of simple Lie algebra. 
\end{prop}

One obtains the following models of the exceptional complex Lie 
algebras, where we denote by $V_{\l}$ the irreducible $\fh$-module 
of highest weight $\l$, and we follow the conventions
of \cite{bou} for the fundamental weights:
\medskip
$$\begin{array}{ccccc}
 \fh & & T & & \fg \\
\fsl_2 & & V_{3\o_1} & & \fg_2 \\ 
\fsp_6 & & V_{\o_3} & & \f4 \\ 
\fsl_6 & & V_{\o_3} & & \fe_6 \\ 
\fso_{12} & & V_{\o_6} & & \fe_7 \\ 
\fe_7 & & V_{\o_1} & & \fe_8 .
\end{array}$$

\subsection{The non-fundamental case}

To complete our proof of the classification of
complex simple Lie algebras, we prove that the only simple Lie algebras for
which the adjoint  representation is not fundamental are $\fsl_n$ and
$\fsp_{2m}$. 
 
Consider   again the adjoint grading $\fg=\fg_{-2}\op\fg_{-1}\op
\fg_{0}\op\fg_{1}\op\fg_{2}$ of the Lie algebra $\fg$. By \cite{LM0},
Proposition 4.2, $\fg_1$ has at most two irreducible components  
as an $\fh$-module, where $\fh\subset\fg_0$
is the semisimple part of $\fg_0$, and the  smallest weights
of these $\fh$-modules are the
simple roots
$\a_i$ such that
$\tilde\a-\a_i$ is still a root.

\medskip\noindent {\bf First case}. Suppose that $\fg_1$ is 
irreducible. If the adjoint representation is not fundamental,
$\tilde\a$ is a multiple of some fundamental weight $\o$
and thus $X=v_d(Z)$ where $Z\subset \ppp V_{\o}$ is the closed orbit and
$d>1$.

We give two proofs of this case:

First proof: $\tbase |II|=\emptyset$ because a Veronese re-embedding of
any variety contains
no lines and a Veronese re-embedding of a homogeneous variety
is such that its ideal is generated  in degree two.
However, we know that the quadratic 
equations of $\tbase |II|$ in $S^2T_1^*$ are given by $\fh$.
Thus $\fh=S^{(2)}T_1^*$, since the other components of $S^2T_1^*$
all vanish on the closed orbit of $\PP T_1$. So $\fh$ itself is 
not fundamental, and by induction on the dimension we conclude that
$\fh$ must be a symplectic Lie algebra and $T_1$ is its natural 
representation.

Second proof:  By \cite{lanci} 3.10,
$  |II|=  S^2T_1^*$.  Since $N_2\simeq \fh \oplus \fg_1\oplus\fg_2$,
we see     $|II|= \fh$ and in particular that $\fh\simeq S^2T^*_1$, i.e,
that
$\fh={\mathfrak sp}_{2m}$ for some $m$ and $T_1$ is its standard 
representation. Since the $\fh$-module structure of $T_1$ 
completely determines $\fg$, this implies that $\fg$ is also
a symplectic Lie algebra. 

\smallskip
In this case our construction works as well, only it is degenerate,
of degree two. Taking $Y=\emptyset$, the rational map $\PP^n
\dashrightarrow \PP^N$ is by the quadrics vanishing on $Y$, that is, 
all quadrics. If one adopts  the convention that $\tau (Y)=\tau
(Y)_{sing}=\PP T_1$, then this construction becomes consistent with the
adjoint algorithm.

\medskip\noindent {\bf Second case}.
Suppose that $\fg_1$ is reducible, hence the sum of two dual 
$\fg$-modules $U$ and $U^*$ (this duality comes from the 
bracket $[\fg_1,\fg_1]\rightarrow \fg_2$).

\begin{lemm} $\PP U$ is $H$-homogeneous. \end{lemm}

\begin{proof} Let the simple root $\a$ be the smallest 
weight of $U$. Then the highest root $\tilde\a$ has coefficient 
one on $\a$, and the weights $\g$ of $U$ are the other positive 
roots with this property. We write $\g=\sum_{\b\;simple}n_{\b}(\g)\b$. 
Note that $n(\a,\tilde\a)=1$ so that
$$ 
  n(\tilde\a,\a)=2+\sum_{\b\ne\a}n_{\b}(\tilde\a)n(\b,\a)>0.
$$
The coefficients $n(\b,\a)$ are the Cartan numbers, which are 
negative for $\b\ne\a$. 
Now $n_{\b}(\g)\le n_{\b}(\tilde\a)$
  for any positive root $\g$ and any
simple root $\b$ so the equation above 
  implies that $n(\tilde\a,\g)>0$ for each root $\g$ that 
is a weight of $U$. Therefore, $\g -k\a$ is a root
for some $k>0$ and in particular $\g-\a$ is   a root. 

Geometrically, letting
$Y$ denote the smallest $H$-orbit in $\ppp U$
and $v_\a\in U$ a vector of weight $\a$,
   this implies that  $\hat T_{[v_{\a}]}Y
=  
\fh v_{\a} =U$, i.e.,  that $Y=\ppp U$.
 \end{proof}

\medskip
Let $v=u+u^*\in\fg_1=U\op U^*$.  Now $v\in \tbase |II|$  if and only
$\o^*(v,Xv)=
2\langle u,Xu^*\rangle=0 $
  for every $X\in \fh$. By the lemma, this is the case if and only
if $u$ or $u^*$ is zero, i.e.,  
$$\tbase |II|=\ppp U\sqcup \ppp U^*\subset\PP (U\op U^*).$$
Moreover,    $I_2(\tbase |II|)=U\ot U^*\subset S^2\fg_1^*=S^2U\op U\ot
U^*\op S^2U^*$.  Thus $U\ot U^*=\fh$. Finally, since $\fg$ is completely
determined  by the $\fh$-module structure of $\fg_1$, we conclude that 
$\fg=\fsl_m$ for some integer $m$. Our proof of the
classification of complex simple Lie algebras is now complete.

\smallskip

If we adopt the convention that $\tau (Y)$ is 
the quartic hypersurface 
  $\langle u,u^*\rangle ^2=0$, which is a double 
quadric,  and  $I_3(\tau(Y)_{sing})$ is the space of derivatives of this
quartic, then
 the adjoint algorithm also works in this case. (However usually
one would take the tangential variety of $Y$ to be $Y$ itself.) 

\smallskip

For example,   for the adjoint variety of $\fsl_3$ (the
projectivization of the space of matrices of order three, rank one 
and trace zero)  we have the following quartic parametrization:

$$(a,b,c,d)\mapsto \begin{pmatrix}   
     a^3d+a^2bc &  -a^2bd-ab^2c &   b^2c^2-a^2d^2 \\  
   2a^3c  &          -2a^2bc    &   2abc^2-2a^2cd \\
         a^4      &      -a^3b   &     a^2bc-a^3d
\end{pmatrix}.   $$

\bigskip

Joseph M. Landsberg 

Mathematics Department,

Georgia Institute of Technology,

Atlanta, GA 30332-0160

UNITED STATES

{\rm E-mail}: jml@math.gatech.edu

\bigskip
Laurent Manivel

Institut Fourier, UMR 5582 du CNRS

Universit\'e Grenoble I

BP 74

38402 Saint Martin d'H\`eres cedex

FRANCE

{\rm E-mail}: Laurent.Manivel@ujf-grenoble.fr 

\end{document}